\documentclass[aps,pre,showpacs,preprint,superscriptaddress]{revtex4-1}
\usepackage{graphicx}
\usepackage{subfigure}
\usepackage{amsmath}
\usepackage{ulem}
\usepackage{bm}
\usepackage{amssymb}
\usepackage{color}
\usepackage{epstopdf}
\usepackage{epsfig,dsfont,multirow}

\begin{document}

\title{Stability range of parameters at fixed points for a class of complex dynamics}
\author{Zhen-Hua Feng}
\affiliation{Department of Physics, Beijing Normal University, Beijing 100875, China}
\author{Hai-Bo Sang}
\affiliation{Department of Physics, Beijing Normal University, Beijing 100875, China}
\author{B. S. Xie \footnote{Corresponding author. Email: bsxie@bnu.edu.cn}}
\affiliation{Key Laboratory of Beam Technology of the Ministry of Education, and College of Nuclear Science and Technology, Beijing Normal University, Beijing 100875, China}
\affiliation{Institute of Radiation Technology, Beijing Academy of Science and Technology, Beijing 100875, China}
\date{\today}
\begin{abstract}
We study the parameters range for the fixed point of a class of complex dynamics with the rational fractional function as $R_{n,a,c}(z)=z^n+\frac{a}{z^n}+c$, where $n=1,2,3,4$ is specified, $a$ and $c$ are two complex parameters. The relationship between two parameters, $a$ and $c$, is obtained at the fixed point.
Moreover the explicit expression of the parameter $a$ and $c$ in terms of $\lambda$ is derived, where $\lambda$ is the derivative function at fixed point. The parameter regimes for the stability of the fixed point are presented numerically for some typical different cases.
\end{abstract}
\maketitle

\section{Introduction}

Rational fractional complex dynamics is the study of iteration, which is the process of repeatedly applying a rational fractional $R(z)$ to itself. For example, we select a value $z_0$ in the domain of $R(z)$ and then use $R(z)$ to continuously iterate at the point $z_0$ to obtain the sequence
\begin{equation}
    \left\{z_0, z_1=R(z_0), z_2=R(z_1), \dots, z_n=R(z_{n-1}),\dots \right\}.
    \label{eq:orbit0}
\end{equation}
Note that the sequence is the orbit of $z_0$ under $R(z)$. $z_0$ is a fixed point if $z_0=z_1=z_2=\dots$. The derivative of $R(z)$ at $z_0$ is the multiplier $\lambda$, i.e. $R^\prime(z_0)=\lambda$. Fixed points can be classified according to $\lambda$ as follows~\cite{ComDyn}:
\begin{itemize}
    \item Superattraction: $\lambda=0$
    \item Attraction: $0<|\lambda|<1$.
    \item Repelling: $|\lambda|>1$.
    \item Rationally neutral: $|\lambda|=1$ and $|\lambda^n|=1$ for some integer $n$.
    \item Irrationally neutral: $|\lambda|=1$ but $|\lambda^n|$ is never 1.
\end{itemize}
These superattraction, attraction and repelling fixed points make up the Fatou set, the complement of Fatou set is Julia set~\cite{AITJF}.

For $R(z)=z^2+c$, those parameters c that can generate the Julia set at the superattraction point $z=0$ form the famous Mandelbrot set~\cite{ComDyn,AITJF,MJFO, MSPPCRF}. Furthermore, the rational fractions $R_{n,a,0}=z^n+\frac{a}{z^n}$ and $R_{n,a,c}=z^n+\frac{a}{z^n}+c$ have been studied in detail~\cite{MSPPCRF,GBBM}. However, the relationship between parameters is not obtained.

We study the rational fractional $R_{n,a,c}=z^n+\frac{a}{z^n}+c$ and obtain the relationship between the parameters $a$, $c$ and $\lambda$, which denoted $F(a,c,\lambda)$. The relationship between the parameters is precise when $n=1,2,3,4$ and approximate when $n>4$. When one parameter is fixed, $F(a,c,\lambda)$ maps the unit circle in the $\lambda$ plane to the stability range of the other parameter.

\section{The relationship of the parameters}

For $R(z)=z^n+\frac{a}{z^n}+c$, the fixed point satisfies the following equation
\begin{equation}
    z=z^n+\frac{a}{z^n}+c,
    \label{eq:ylfs-1}
\end{equation}
and the multiplier $\lambda$ satisfies
\begin{equation}
    \frac{\lambda}{n}z=z^n-\frac{a}{z^n}.
    \label{eq:ylfs0}
\end{equation}
Substituting Eq.\eqref{eq:ylfs-1} into Eq.\eqref{eq:ylfs0}, we obtain
\begin{equation}
    \frac{\lambda}{n}\left( z^n+\frac{a}{z^n}+c \right)=z^n-\frac{a}{z^n}.
   \label{eq:ylfs1}
\end{equation}
After a simple calculation, we obtain
\begin{equation}
    \left( 1-\frac{\lambda}{n} \right)z^{2n}-\frac{\lambda}{n}cz^{n}-\left( 1+\frac{\lambda}{n} \right)a=0.
   \label{eq:ylfs2}
\end{equation}
The Eq. \eqref{eq:ylfs2} is easily solved for $z^n$ as that
\begin{equation}
    z^{n}=\frac{c}{2\left( 1-\frac{\lambda}{n} \right)}\left[ \frac{\lambda}{n}\pm\sqrt{\frac{\lambda^2}{n^2}+\frac{4a}{c^2}\left( 1-\frac{\lambda^2}{n^2} \right)}\right].
   \label{eq:ylfs3}
\end{equation}
Now we define
\begin{equation}
    \xi=\frac{c}{1-\frac{\lambda^2}{n^2}}\left[ 1\pm\sqrt{\frac{\lambda^2}{n^2}+\frac{4a}{c^2}\left( 1-\frac{\lambda^2}{n^2} \right)} \right],
   \label{eq:ylfs5}
\end{equation}
and combining Eq.\eqref{eq:ylfs-1} and Eq.\eqref{eq:ylfs3}, it yields
\begin{equation}
    2\xi^n-\left( 1+\frac{\lambda}{n} \right)\xi+c=0.
   \label{eq:ylfs6}
\end{equation}

From Eqs.\eqref{eq:ylfs5} and \eqref{eq:ylfs6} one can obtain $c$ easily,
\begin{equation}
    c=\left( 1-\frac{\lambda}{n} \right)\xi-\frac{2a}{\xi^n}.
    \label{eq:ylfs7a}
\end{equation}
Moreover, combining Eqs.\eqref{eq:ylfs6} and \eqref{eq:ylfs7a}, we obtain
    \begin{equation}
        \xi^{2n}-\frac{\lambda}{n}\xi^{n+1}-a=0.
        \label{eq:ylfs7b}
    \end{equation}
When Eq.\eqref{eq:ylfs6} can be solved for $\xi$, then substituting $\xi$ into Eq.\eqref{eq:ylfs7b}, therefore, one concludes that $a$ can be expressed by $c$ and $\lambda$. Similarly, by inserting the $\xi$ as roots of Eq.\eqref{eq:ylfs7b} into Eq.\eqref{eq:ylfs6}, then $c$ would be expressed by $a$ and $\lambda$.

\section{Solvable cases}

If $a=0$ and $c=0$, the rational fraction we studied is meaningless. Therefore, this paper would study three cases of $a=0$ and $c\neq0$, $a\neq0$ and $c=0$, and $a\neq0$ and $c\neq0$. These cases that can be solved are as follows.

\subsection{The case $a=0$ and $c\neq0$}

From Eq.\eqref{eq:ylfs5}, we obtain $\xi\neq 0$ for $a=0$ or $c=0$. if $a=0$ and $c\neq0$, Eq.\eqref{eq:ylfs7b} can be solved for $\xi$, that is
\begin{equation}
    \xi_k=\left( \frac{\lambda}{n} \right)^{\frac{1}{n-1}}e^{i\frac{2k\pi}{n-1}}
    \qquad
    \left( k= 0,1,\dots , n-2 \right).
    \label{eq:ylfs8}
\end{equation}
Combining Eq.\eqref{eq:ylfs8} and Eq.\eqref{eq:ylfs6}, we obtain
\begin{equation}
    c_k=\left( 1-\frac{\lambda}{n} \right)\xi_k
    \qquad
    \left( k= 0,1,\dots , n-2 \right).
    \label{eq:ylfs8a}
\end{equation}

\subsection{The case $a\neq0$ and $c=0$}
If $c=0$ and $a\neq0$,  we solved Eq.\eqref{eq:ylfs6} to obtain
\begin{equation}
    \xi_k=\left( \frac{1+\frac{\lambda}{n}}{2} \right)^{\frac{1}{n-1}}e^{i\frac{2k\pi}{n-1}}
    \qquad
    \left( k= 0,1,\dots , n-2 \right).
    \label{eq:ylfs8b}
\end{equation}
Combining Eq.\eqref{eq:ylfs8b} and Eq.\eqref{eq:ylfs7b}, then parameter $a$ is obtained as
\begin{equation}
    a_k=\frac{1}{2}\left( 1-\frac{\lambda}{n} \right)\xi_k^{n+1}
    \qquad
    \left( k= 0,1,\dots , n-2 \right).
    \label{eq:ylfs8c}
\end{equation}

\subsection{The case $a\neq0$ and $c\neq0$}

It is well known that the rational equations of less than $5$ degrees can be solved. Therefore, this paper obtains the parameter relationship when $n = 1,2,3,4$.

\subsubsection{Parameter relationship when $n=1$}

Let $n=1$ in Eq.\eqref{eq:ylfs6}, then we obtain
\begin{equation}
    2\xi-\left( 1+\lambda \right)\xi+c=0.
   \label{eq:ylfs11}
\end{equation}
Solving it for $\xi$, one has
\begin{equation}
    \xi=-\frac{c}{1-\lambda}.
   \label{eq:ylfs12}
\end{equation}
So from Eq.\eqref{eq:ylfs7b}, it yields
\begin{equation}
    a=\frac{c^2}{1-\lambda}.
   \label{eq:ylfs13}
\end{equation}
And \emph{vice versa}, it is also expressed for $c$ with $a$ as
\begin{equation}
    c=\pm\sqrt{a(1-\lambda)}.
    \label{eq:ylfs13a}
\end{equation}

\subsubsection{Parameter relationship when $n=2$}

In Eq.\eqref{eq:ylfs6}, we let $n=2$, then we obtain
\begin{equation}
    2\xi^2-\left( 1+\frac{\lambda}{2} \right)\xi+c=0.
   \label{eq:ylfs14}
\end{equation}
Solving it for $\xi$, we obtain
\begin{equation}
   \begin{cases}
      \xi_0=\frac{1}{4}\left[ \left( 1+\frac{\lambda}{2} \right)+\sqrt{\left( 1+\frac{\lambda}{2} \right)^2-8c} \right]\\
      \xi_1=\frac{1}{4}\left[ \left( 1+\frac{\lambda}{2} \right)-\sqrt{\left( 1+\frac{\lambda}{2} \right)^2-8c} \right]
  \end{cases}.
   \label{eq:ylfs15}
\end{equation}
Combining Eqs.\eqref{eq:ylfs7b} and \eqref{eq:ylfs15}, it yields
\begin{equation}
    a_k=\xi_k^{4}-\frac{\lambda}{2}\xi_k^{3}
   \qquad
   \left( k=0,1 \right).
   \label{eq:ylfs16}
\end{equation}

Moreover, let $n=2$ in Eq.\eqref{eq:ylfs7b}, then we have
    \begin{equation}
        \xi^{4}-\frac{\lambda}{2}\xi^{3}-a=0.
        \label{eq:ylfs16a}
    \end{equation}

The standard solution of Eq.\eqref{eq:ylfs16a} for $\xi$ can be expressed as:
    \begin{equation}
        \begin{cases}
            \xi_0 = \frac{1}{2} \left[ \left( \frac{\lambda}{4}+\sqrt{2\alpha+\frac{\lambda^2}{16}} \right)+\sqrt{\left( \frac{\lambda}{4}+\sqrt{2\alpha+\frac{\lambda^2}{16}} \right)^2-4\alpha\left( 1+\frac{\lambda}{4\sqrt{2\alpha+\frac{\lambda^2}{16}}} \right)} \right]\\
            \xi_1 = \frac{1}{2} \left[ \left( \frac{\lambda}{4}+\sqrt{2\alpha+\frac{\lambda^2}{16}} \right)-\sqrt{\left( \frac{\lambda}{4}+\sqrt{2\alpha+\frac{\lambda^2}{16}} \right)^2-4\alpha\left( 1+\frac{\lambda}{4\sqrt{2\alpha+\frac{\lambda^2}{16}}} \right)} \right]\\
            \xi_2 = \frac{1}{2} \left[ \left( \frac{\lambda}{4}-\sqrt{2\alpha+\frac{\lambda^2}{16}} \right)+\sqrt{\left( \frac{\lambda}{4}-\sqrt{2\alpha+\frac{\lambda^2}{16}} \right)^2-4\alpha\left( 1-\frac{\lambda}{4\sqrt{2\alpha+\frac{\lambda^2}{16}}} \right)} \right]\\
            \xi_3 = \frac{1}{2} \left[ \left( \frac{\lambda}{4}-\sqrt{2\alpha+\frac{\lambda^2}{16}} \right)-\sqrt{\left( \frac{\lambda}{4}-\sqrt{2\alpha+\frac{\lambda^2}{16}} \right)^2-4\alpha\left( 1-\frac{\lambda}{4\sqrt{2\alpha+\frac{\lambda^2}{16}}} \right)} \right]
        \end{cases},
        \label{eq:ylfs16c}
    \end{equation}
where the definition of $\alpha$ is denoted as
    \begin{equation}
        \alpha = \left[ -\left( \frac{\lambda}{8} \right)^2a+\sqrt{\left( \frac{\lambda}{8} \right)^4a^2+\left( \frac{a}{3} \right)^{3}} \right]^{\frac{1}{3}}
        +
        \left[ -\left( \frac{\lambda}{8} \right)^2a-\sqrt{\left( \frac{\lambda}{8} \right)^4a^2+\left( \frac{a}{3} \right)^{3}} \right]^{\frac{1}{3}}.
        \label{eq:ylfs16b}
    \end{equation}

Using Eq.\eqref{eq:ylfs6} with $\xi$ of Eq.\eqref{eq:ylfs16c}, finally, we obtain the required expression of $c$ based on the $a$ and $\lambda$ as
\begin{equation}
    c_k=\left( 1+\frac{\lambda}{2} \right)\xi_k-2\xi_k^{2}
    \qquad
    (k=0,1,2,3).
    \label{eq:ylfs16d}
\end{equation}

\subsubsection{Parameter relationship when $n=3$}

Similarly, let $n=3$, then Eq.\eqref{eq:ylfs6} becomes
\begin{equation}
    2\xi^3-\left( 1+\frac{\lambda}{3} \right)\xi+c=0.
   \label{eq:ylfs17}
\end{equation}
Through two auxiliary quantities
\begin{equation}
    u_0=\left\{\frac{1}{2}\left[- \frac{c}{2}+\sqrt{\frac{c^2}{4}-4\left( \frac{1+\frac{\lambda}{3}}{6} \right)^3} \right]\right\}^{\frac{1}{3}},
   \label{eq:ylfs23}
\end{equation}
and
\begin{equation}
    v_0=\left\{\frac{1}{2}\left[- \frac{c}{2}-\sqrt{\frac{c^2}{4}-4\left( \frac{1+\frac{\lambda}{3}}{6} \right)^3} \right]\right\}^{\frac{1}{3}},
   \label{eq:ylfs24}
\end{equation}
we can obtain the $\xi$ for the solution of Eq.\eqref{eq:ylfs17} as
\begin{equation}
   \begin{cases}
      \xi_0=u_0+v_0\\
      \xi_1=u_0e^{i\frac{2\pi}{3}}+v_0e^{-i\frac{2\pi}{3}}\\
      \xi_2=u_0e^{i\frac{4\pi}{3}}+v_0e^{-i\frac{4\pi}{3}}
  \end{cases}.
   \label{eq:ylfs25}
\end{equation}
Thus Eqs.\eqref{eq:ylfs7b} and \eqref{eq:ylfs25} yield the required relation
\begin{equation}
    a_k=\xi_k^{6}-\frac{\lambda}{3}\xi_k^{4}
   \qquad
   \left( k=0,1,2 \right).
   \label{eq:ylfs26}
\end{equation}

Moreover, let $n=3$ in Eq.\eqref{eq:ylfs7b}, then we have
\begin{equation}
    \xi^6-\frac{\lambda}{3}\xi^4-a=0.
    \label{eq:ylfs26a}
\end{equation}
Through two auxiliary quantities
\begin{equation}
    p=\left\{ \frac{a}{2}+\left( \frac{\lambda}{9} \right)^3+\sqrt{\frac{a}{2}\left[ \frac{a}{2}+2\left( \frac{\lambda}{9} \right)^3 \right]} \right\}^{\frac{1}{3}},
    \label{eq:ylfs26b}
\end{equation}
and
\begin{equation}
    q=\left\{\frac{a}{2}+\left( \frac{\lambda}{9} \right)^3-\sqrt{\frac{a}{2}\left[ \frac{a}{2}+2\left( \frac{\lambda}{9} \right)^3 \right]} \right\}^{\frac{1}{3}},
    \label{eq:ylfs26c}
\end{equation}
we can obtain the $\xi$ for the solution of Eq.\eqref{eq:ylfs26a} as
\begin{equation}
    \begin{cases}
        \xi_0=\sqrt{\frac{\lambda}{9}+p+q}\\
        \xi_1=-\sqrt{\frac{\lambda}{9}+p+q}\\
        \xi_2=\sqrt{\frac{\lambda}{9}+pe^{i\frac{2\pi}{3}}+qe^{-i\frac{2\pi}{3}}}\\
        \xi_3=-\sqrt{\frac{\lambda}{9}+pe^{i\frac{2\pi}{3}}+qe^{-i\frac{2\pi}{3}}}\\
        \xi_4=\sqrt{\frac{\lambda}{9}+pe^{i\frac{4\pi}{3}}+qe^{-i\frac{4\pi}{3}}}\\
        \xi_5=-\sqrt{\frac{\lambda}{9}+pe^{i\frac{4\pi}{3}}+qe^{-i\frac{4\pi}{3}}}
    \end{cases}.
    \label{eq:ylfs26d}
\end{equation}
Combining Eqs.\eqref{eq:ylfs26d} and \eqref{eq:ylfs7a}, finally, it yields
\begin{equation}
    c_k=\left( 1+\frac{\lambda}{3} \right)\xi_k-2\xi_k^{3}
    \qquad
    \left( k=0,1,2,3,4,5 \right)
    .
    \label{eq:ylfs26e}
\end{equation}

\subsubsection{Parameter relationship when $n=4$}

Let us to look at the case of $n=4$, now Eq. \eqref{eq:ylfs6} is
\begin{equation}
    2\xi^4-\left( 1+\frac{\lambda}{4} \right)\xi+c=0.
   \label{eq:ylfs27}
\end{equation}
Define an auxiliary
\begin{align}
   \alpha&=\left[\left( \frac{1+\frac{\lambda}{4}}{8} \right)^2+\sqrt{\left( \frac{1+\frac{\lambda}{4}}{8} \right)^2-\left( \frac{c}{6} \right)^3}\right]^{\frac{1}{3}}\notag\\
   &+\left[\left( \frac{1+\frac{\lambda}{4}}{8} \right)^2-\sqrt{\left( \frac{1+\frac{\lambda}{4}}{8} \right)^2-\left( \frac{c}{6} \right)^3}\right]^{\frac{1}{3}},
   \label{eq:ylfs29ans}
\end{align}
to solve the Eq.\eqref{eq:ylfs27} for $\xi$, one obtain the following solutions as
\begin{equation}
   \begin{cases}
       \xi_0=\frac{1}{2}\left[ \sqrt{2\alpha}+\sqrt{-2\alpha+\frac{1+\frac{\lambda}{4}}{\sqrt{2\alpha}} }\right]\\
   \xi_1=\frac{1}{2}\left[ \sqrt{2\alpha}-\sqrt{-2\alpha+\frac{1+\frac{\lambda}{4}}{\sqrt{2\alpha}} }\right]\\
   \xi_2=\frac{1}{2}\left[ -\sqrt{2\alpha}+\sqrt{-2\alpha-\frac{1+\frac{\lambda}{4}}{\sqrt{2\alpha}} }\right]\\
   \xi_4=\frac{1}{2}\left[ -\sqrt{2\alpha}-\sqrt{-2\alpha-\frac{1+\frac{\lambda}{4}}{\sqrt{2\alpha}} }\right]
   \end{cases}.
   \label{eq:ylfs33}
\end{equation}
Similarly, through the $\xi$ in Eqs.\eqref{eq:ylfs7b} and \eqref{eq:ylfs33}, finally, it yields
\begin{equation}
    a_k=\xi_k^{8}-\frac{\lambda}{4}\xi_k^{5}
   \qquad
   \left( k=0,1,2,3 \right).
   \label{eq:ylfs34}
\end{equation}

\section{Approximated cases}

When $n$ is large, the equations in the form of Eq.\eqref{eq:ylfs6} or/and Eq.\eqref{eq:ylfs7b} can not be solved for $\xi$. In a certain condition, however, they do exist some approximate solutions, which are helpful to the problem of finding the stability parameters ranges of the fixed points in present study.

\subsection{An approximate expression for parameter $a$}

When $n\geq 5$ and meanwhile if $|c|\gg 2^{\frac{1}{1-n}}|1+\frac{\lambda}{n}|^{\frac{n}{n-1}}$, then Eq.\eqref{eq:ylfs6} can be approximately solved for $\xi$ to lead
\begin{equation}
    \xi_k=\left( \frac{c}{2} \right)^{i\frac{(2k-1)\pi}{n}}
   \qquad
   \left( k=1,2,\cdots , n \right).
    \label{eq:ylfs6a}
\end{equation}
With $\xi$ in Eq.\eqref{eq:ylfs7b}, we have
\begin{equation}
    a_k=\xi_k^{2n}-\frac{\lambda}{n}\xi_k^{n+1}
   \qquad
   \left( k=1,2,\cdots , n \right).
    \label{eq:ylfs6b}
\end{equation}

On the other hand, if $0<|c|\ll 2^{\frac{1}{1-n}}|1+\frac{\lambda}{n}|^{\frac{n}{n-1}}$, then Eq.\eqref{eq:ylfs6} can be approximately solved for $\xi$ to obtain
\begin{equation}
    \xi=\frac{c}{1+\frac{\lambda}{n}},
    \label{eq:ylfs6c}
\end{equation}
Substituting Eq.\eqref{eq:ylfs6c} into Eq.\eqref{eq:ylfs7b}, we obtain
\begin{equation}
    a=\xi^{2n}-\frac{\lambda}{n}\xi^{n+1}.
    \label{eq:ylfs6d}
\end{equation}

\subsection{An approximate expression for parameter $c$}

Similarly, when $n\geq 4$ and meanwhile if $|a|\gg|\frac{\lambda}{n}|^{\frac{2n}{2n-1}}$, Eq.\eqref{eq:ylfs7b} can be also solved approximately,  that is
\begin{equation}
    \xi_k=a^{\frac{1}{2n}}e^{i\frac{k\pi}{n}}
   \qquad
   \left( k=1,2,\cdots , 2n \right).
    \label{eq:ylfs7c}
\end{equation}
Substituting these $\xi$ into Eq.\eqref{eq:ylfs6} would yield the following result
\begin{equation}
    c=(1+\frac{\lambda}{n})\xi_k-2\xi_k^{n}
   \qquad
   \left( k=1,2,\cdots , 2n \right).
    \label{eq:ylfs7d}
\end{equation}

Similarly, if $0<|a|\ll|\frac{\lambda}{n}|^{\frac{2n}{2n-1}}$, Eq.\eqref{eq:ylfs7b} can be also solved approximately, thus $\xi$ is
\begin{equation}
    \xi=\left(-\frac{na}{\lambda} \right)^{\frac{1}{n+1}}.
    \label{eq:ylfs7e}
\end{equation}
Substituting the $\xi$ into Eq.\eqref{eq:ylfs6} would lead to $c$ as
\begin{equation}
    c=\left( 1+\frac{\lambda}{n} \right)\xi-2\xi^n.
    \label{eq:ylfs7f}
\end{equation}

\section{Numerical results for parameter ranges}

In the previous section we obtained the correlation between parameters at fixed points, and in this section we would plot the range of parameters in several special cases based on the analytical results mentioned above. Note that we omit the stability range of the parameter when it can take infinity. Since $R_{n,a,c}(z)=z^n+\frac{a}{z^n}+c$ has two parameters, If one parameter is being fixed, then $F(a,c,\lambda)$ can map the unit circle in the $\lambda$ plane to the stability range of the other parameter.

\subsection{The case $a=0$ or $c=0$}

\begin{figure}[!ht]
    \begin{center}
        \includegraphics[scale=0.5]{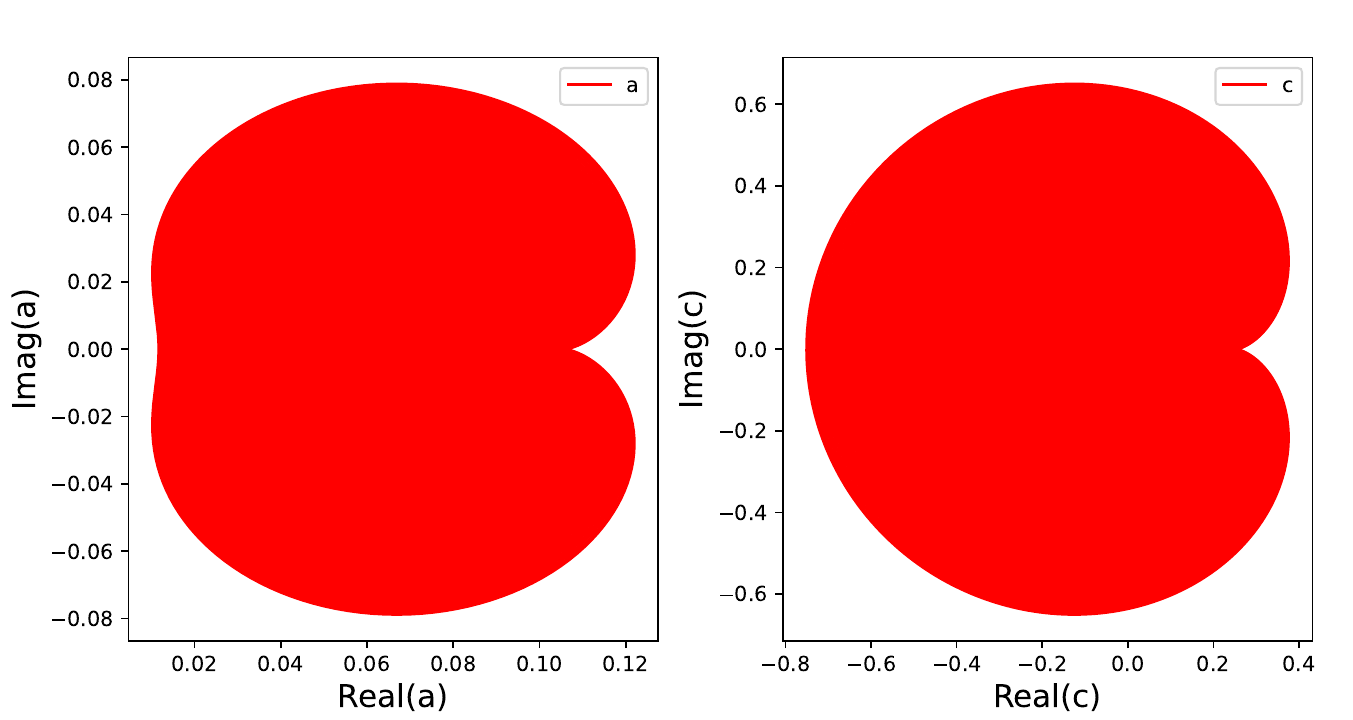}
    \end{center}
    \caption{The stability range of parameters $a$ and $c$ in the case of $n=2$, where $c=0$ is fixed for the range of $a$ (left panel) and $a=0$ is fixed for the range of $c$ (right panel).}
    \label{fig:aorc0}
\end{figure}

The stability range of parameters $a$ or $c$ are obtained for any integer $n$ when $c=0$ or $a=0$. Let $\lambda$ takes the value within the unit circle and $c=0$, then the stability range of $a$ is obtained. Similarly, the range of $c$ can also be obtained for $a=0$. shown in fig.\ref{fig:aorc0} is the case of $n=2$, in which it is comparable and confirmable by the results of appearing the main part of the Mandelbrot set for the $c$. This also proves the validity of our theoretical treatment in this study.

\subsection{The case $a\neq0$ or $c\neq0$}

\subsubsection{Range of parameters when $n=1$}

\begin{figure}[!ht]
    \begin{center}
        \includegraphics[scale=0.5]{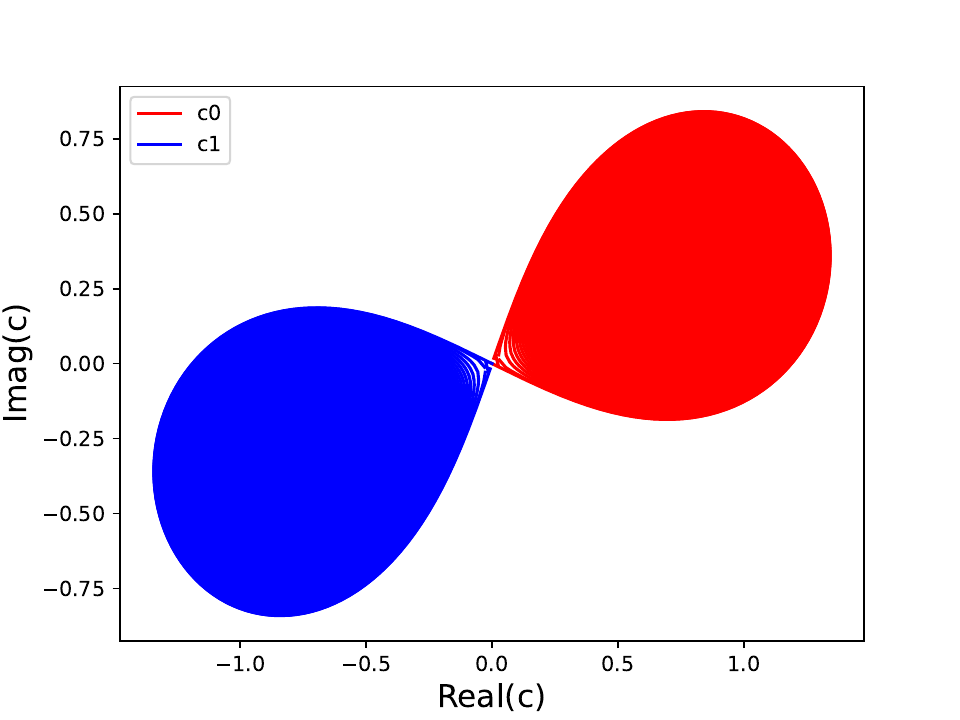}
    \end{center}
    \caption{The stability range of parameter $c$ in the case of $n = 1$ when $a = \frac{\sqrt{2}}{2}(1+i)$ is fixed.}
    \label{fig:n1}
\end{figure}

The stability range of parameter $c$ is finite when $a$ is finite. Let $\lambda$ takes the value within the unit circle, then the stability range of $c$ is shown in Fig.\ref{fig:n1}. By the way it is easily to look at that $c_0=\sqrt{a(1-\lambda)}$ and $c_1=-\sqrt{a(1-\lambda)}$.

\subsubsection{Range of parameters when $n=2$}

\begin{figure}[!ht]
    \begin{center}
        \includegraphics[scale=0.5]{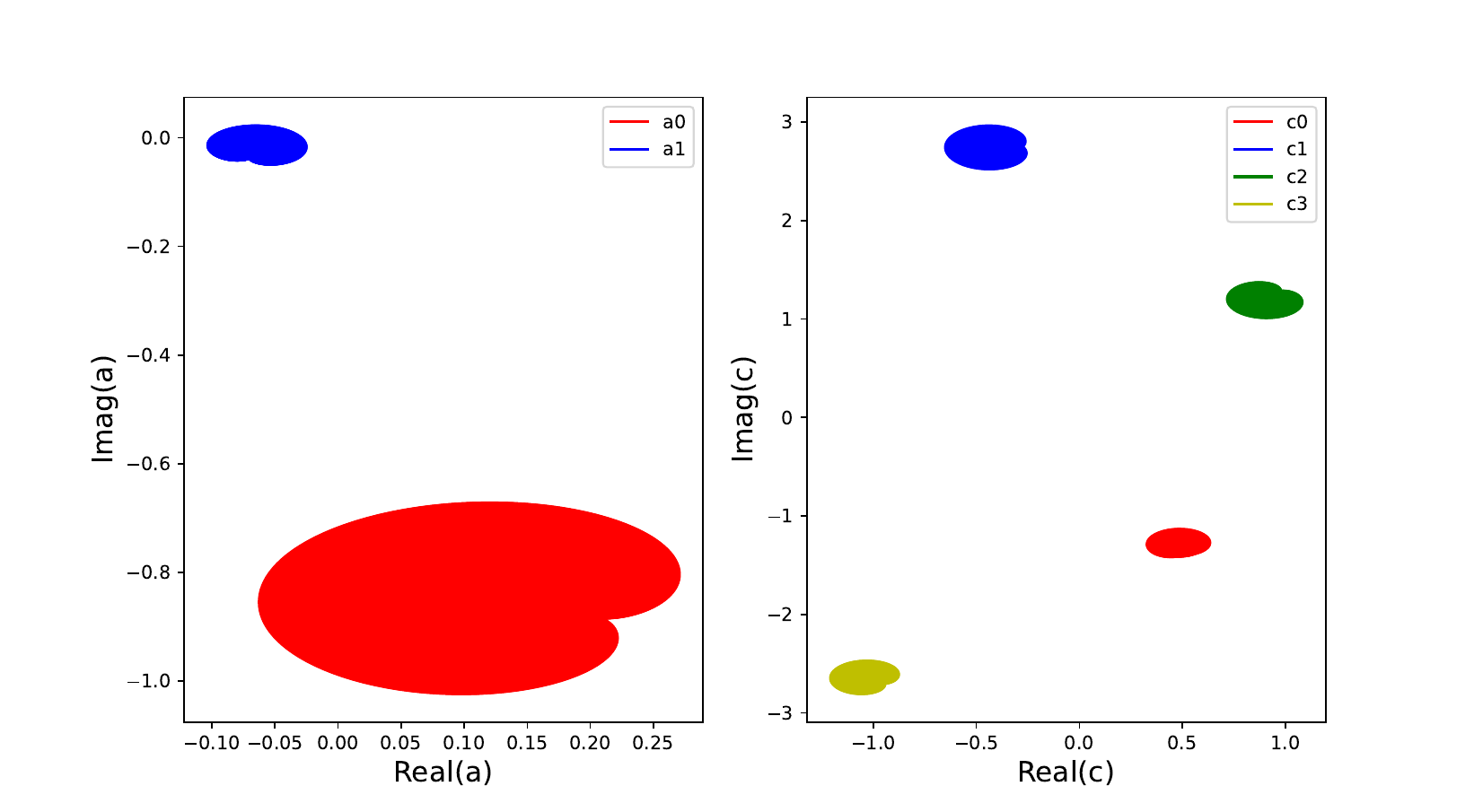}
    \end{center}
    \caption{The stability range of parameters $a$ and $c$ in the case of $n=2$. $c=-\frac{1}{2}+i\frac{\sqrt{3}}{2}$ is fixed for the range of $a$ (left panel) and $a=-\frac{\sqrt{2}}{2}+i\frac{\sqrt{2}}{2}$ is fixed for the range of $c$ (right panel).}
    \label{fig:n2}
\end{figure}

Figure \ref{fig:n2} gives the range of parameters $a$ and $c$, when $n=2$. We fix $c=-\frac{1}{2}+i\frac{\sqrt{3}}{2}$ and $\lambda$ take the value in the unit circle to obtain the range of $a$. Two ranges of $a$ corresponds to two sets of solutions of $\xi$ in Eq.\eqref{eq:ylfs16}. For the range of $c$, we fix $a=-\frac{\sqrt{2}}{2}+i\frac{\sqrt{2}}{2}$ and four sets of $\xi$ are used in Eq.\eqref{eq:ylfs16d}.

\subsubsection{Range of parameters when $n=3$}

\begin{figure}[!ht]
    \begin{center}
        \includegraphics[scale=0.5]{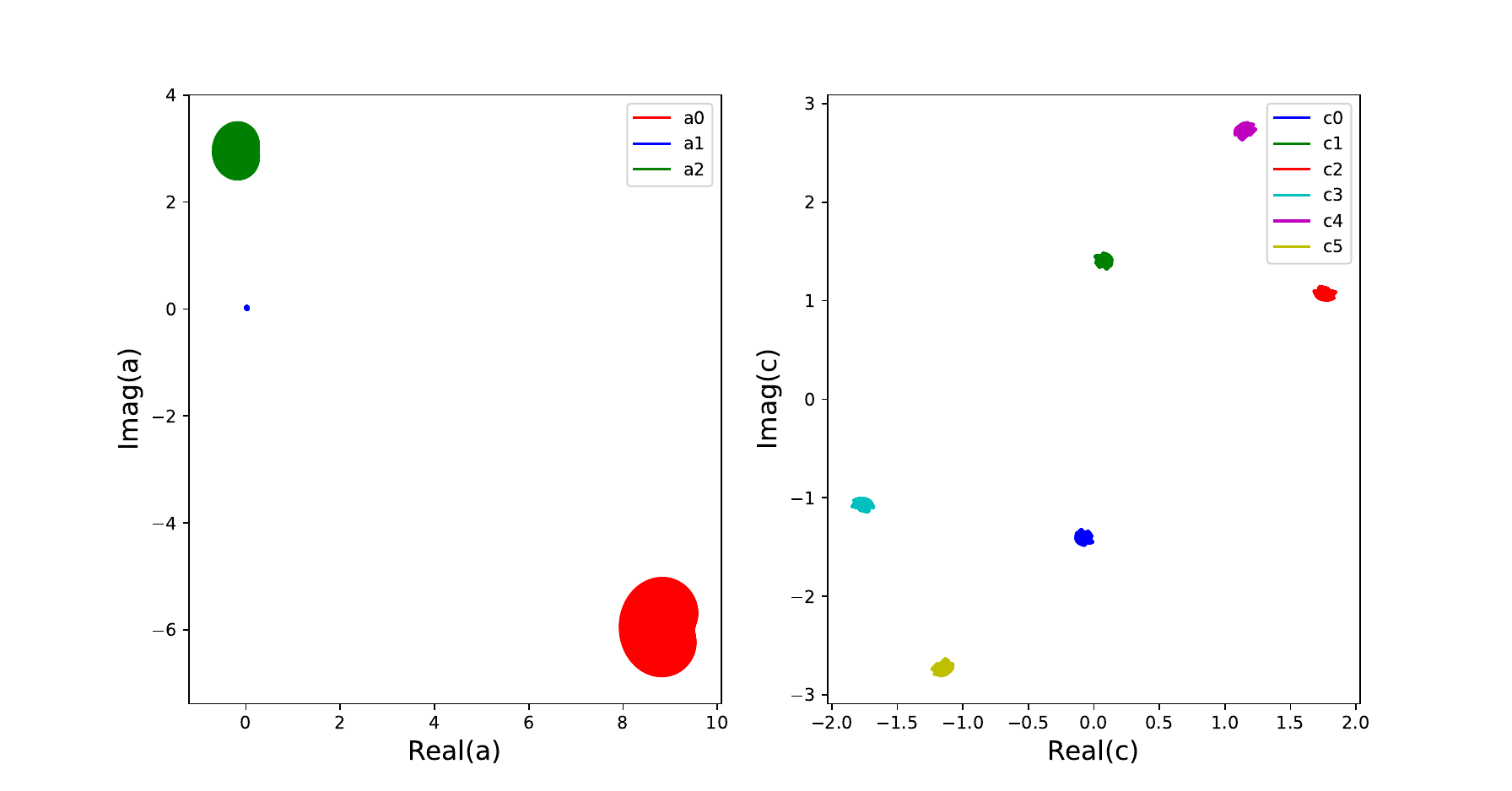}
    \end{center}
\caption{The stability range of parameter $a$(left panel) in the case of $n=3$ when $c=-\frac{\sqrt{2}}{2}+i\frac{\sqrt{2}}{2}$ is fixed and $a=-\frac{1}{2}+i\frac{\sqrt{3}}{2}$ is fixed for the range of $c$(right panel).}
    \label{fig:n3}
\end{figure}

There is only one equation \eqref{eq:ylfs26}, when $n=3$. Let $c=-\frac{\sqrt{2}}{2}+i\frac{\sqrt{2}}{2}$, then we can obtain the stability ranges of $a$, which composed by $a0$, $a1$ and $a3$. Let $a=-\frac{1}{2}+i\frac{\sqrt{3}}{2}$ and $\lambda$ take the value in the unit circle , then the range of $c$ is obtained.

\subsubsection{Range of parameters when $n=4$}

\begin{figure}[!ht]
    \begin{center}
        \includegraphics[scale=0.5]{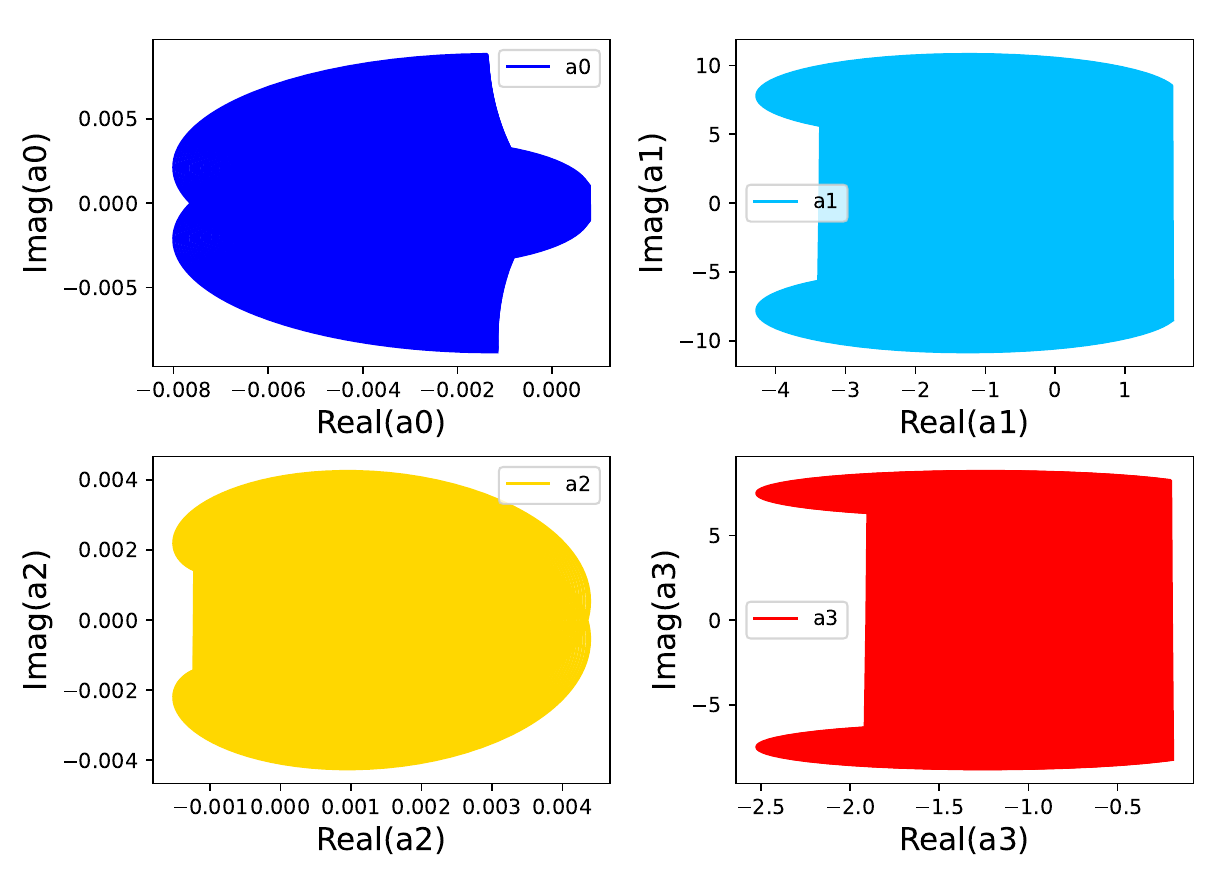}
    \end{center}
    \caption{The stability range of parameters $a$ in the case of $n=4$ when $c=\frac{\sqrt{2}}{2}+i\frac{\sqrt{2}}{2}$ is fixed.}
    \label{fig:n4}
\end{figure}

When $n=4$, there is only equation \eqref{eq:ylfs26} too. Let $c=\frac{\sqrt{2}}{2}+i\frac{\sqrt{2}}{2}$ and $\lambda$ take the value in the unit circle, then the range of $a$ is obtained, which composed by $a0$, $a1$, $a2$ and $a3$.

\subsection{The approximate range of parameters when $n=5$}

\begin{figure}[!ht]
    \begin{center}
        \includegraphics[scale=0.5]{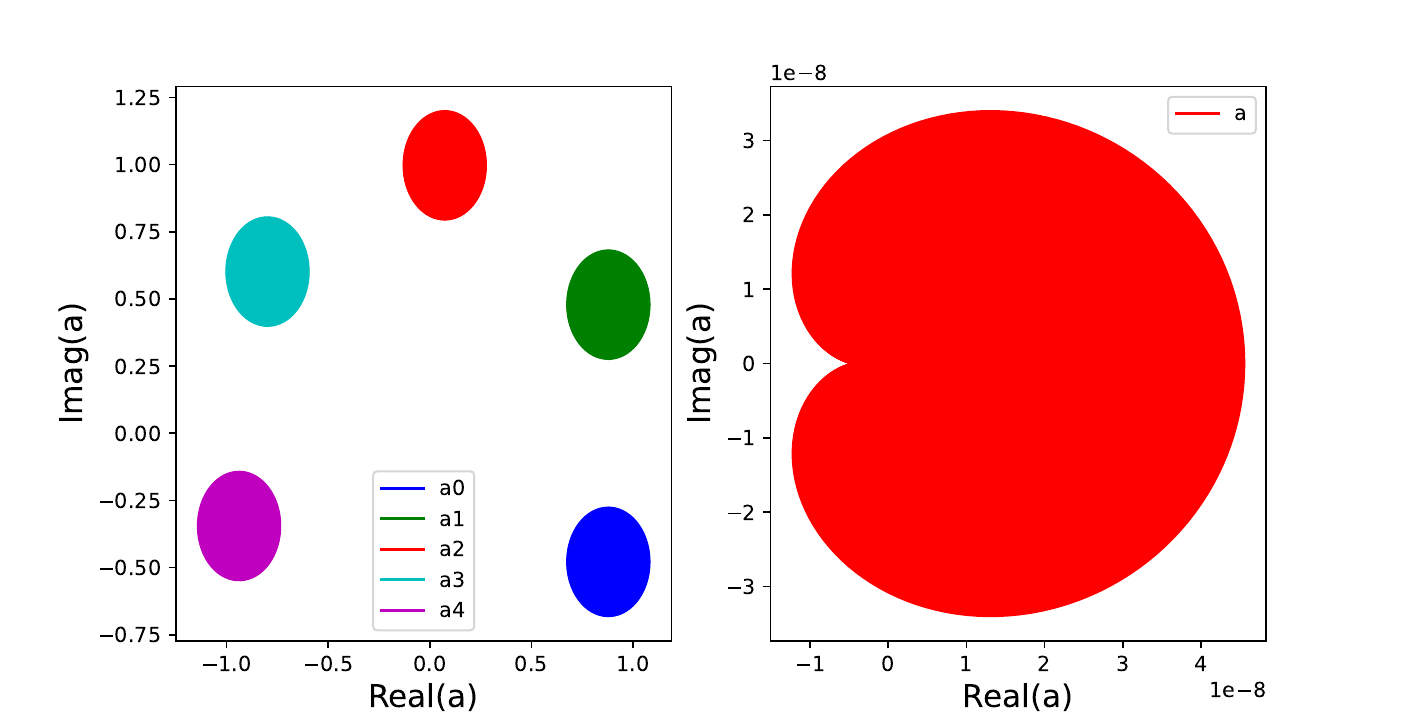}
    \end{center}
    \caption{The approximate stability range of parameters $a$ in the case of $n=5$ when $c=16$ is fixed for the left panel and $c=\frac{1}{16}$ for the right panel.  }
    \label{fig:ap1}
\end{figure}

We cannot get an exact solution when $n=5$, but with the help of an approximate solution, we can roughly obtain the stability range of $a$. Let $\lambda$ take the value in the unit circle, then the ranges of $a$ are obtained for $c=16$ and $c=\frac{1}{16}$.

\section{Conclusion}

In this paper, we obtain the exact or/and approximate analytical expressions about the stability parameter range of the fixed point for $R_{n,a,c}=z^n+\frac{a}{z^n}+c$ when the integer $n$ is small or/and large, respectively. Using these important analytical formulas, we can determine the stability range of the parameters conveniently via numerical calculations, in particularly, for the cases when $n=1,2,3$ and $4$.

\appendix

\section{Derivation of parameters $a$ and $c$ with respect to auxiliary quantity $\xi$}\label{append}

In the previous section we obtained the Eq.\eqref{eq:ylfs3}, that is
\begin{equation}
    z^{n}=\frac{c}{2\left( 1-\frac{\lambda}{n} \right)}\left[ \frac{\lambda}{n}\pm\sqrt{\frac{\lambda^2}{n^2}+\frac{4a}{c^2}\left( 1-\frac{\lambda^2}{n^2} \right)}\right].
   \label{eq:ylfsApp0}
\end{equation}
The product of the reciprocal of the Eq.\eqref{eq:ylfsApp0} and $a$ is
\begin{equation}
    \frac{a}{z^n}=\frac{c}{2\left( 1+\frac{\lambda}{n} \right)}\left[ -\frac{\lambda}{n}\pm\sqrt{\frac{\lambda^2}{n^2}+\frac{4a}{c^2}\left( 1-\frac{\lambda^2}{n^2} \right)} \right].
   \label{eq:ylfsApp1}
\end{equation}
Combining Eqs.\eqref{eq:ylfsApp0} and \eqref{eq:ylfsApp1}, it yields
\begin{equation}
    z^n+\frac{a}{z^n}+c
    =\frac{c}{1-\frac{\lambda^2}{n^2}}\left[ 1\pm\sqrt{\frac{\lambda^2}{n^2}+\frac{4a}{c^2}\left( 1-\frac{\lambda^2}{n^2} \right)} \right].
    \label{eq:ylfsApp2}
\end{equation}
Thereby Eq.\eqref{eq:ylfsApp0} can be rewritten as
\begin{equation}
    z^n=-\frac{c}{2}+\left( 1+\frac{\lambda}{n} \right)\frac{c}{2\left( 1-\frac{\lambda^2}{n^2} \right)}\left[ 1\pm\sqrt{\frac{\lambda^2}{n^2}+\frac{4a}{c^2}\left( 1-\frac{\lambda^2}{n^2} \right)} \right].
    \label{eq:ylfsApp3}
\end{equation}
On the other hand, we note that the $n$-th power of Eq.\eqref{eq:ylfs-1} reads
\begin{equation}
    z^n=\left( z^n+\frac{a}{z^n}+c \right)^{n}.
    \label{eq:ylfsApp4}
\end{equation}
So if we define $\xi$ as
\begin{equation}
    \xi=\frac{c}{1-\frac{\lambda^2}{n^2}}\left[ 1\pm\sqrt{\frac{\lambda^2}{n^2}+\frac{4a}{c^2}\left( 1-\frac{\lambda^2}{n^2} \right)} \right],
    \label{eq:ylfsApp5}
\end{equation}
and substituting Eqs.\eqref{eq:ylfsApp2}, \eqref{eq:ylfsApp3} and \eqref{eq:ylfsApp5} into Eq.\eqref{eq:ylfsApp4}, then, the Eq.\eqref{eq:ylfs6} can be obtained readily as
\begin{equation}
    2\xi^n-\left( 1+\frac{\lambda}{n} \right)\xi+c=0.
   \label{eq:ylfsApp6}
\end{equation}

Moreove, after a simple calculation to Eq.\eqref{eq:ylfsApp5}, it yields
\begin{equation}
    \left( 1-\frac{\lambda^2}{n^2} \right)\left[ \xi\left( 1+\frac{\lambda}{n} \right)-c \right]\left[ \xi\left( 1-\frac{\lambda}{n} \right)-c \right]
    =4a\left( 1-\frac{\lambda^2}{n^2} \right),
    \label{eq:ylfsApp7}
\end{equation}
and from this, it is not difficult to get
\begin{equation}
    a=\frac{1}{4}\left[ \xi\left( 1+\frac{\lambda}{n} \right)-c \right]\left[ \xi\left( 1-\frac{\lambda}{n} \right)-c \right].
    \label{eq:ylfsApp8}
\end{equation}
Furthermore, by combining Eqs.\eqref{eq:ylfsApp6} and \eqref{eq:ylfsApp8}, we have
\begin{equation}
    a=\frac{1}{2}\xi^n\left[ \left( 1-\frac{\lambda}{n} \right)\xi-c \right],
    \label{eq:ylfsApp9}
\end{equation}
which is an expression of parameter $a$ with respect to parameter $c$ and auxiliary quantity $\xi$. Now, as a \emph{vice versa}, the expression of parameter $c$ with respect to parameter $a$ and auxiliary quantity $\xi$ can be obtained as
\begin{equation}
    c=\left( 1-\frac{\lambda}{n} \right)\xi-\frac{2a}{\xi^n}.
    \label{eq:ylfsApp10}
\end{equation}
Finally, by substituting Eq.\eqref{eq:ylfsApp10} into Eq.\eqref{eq:ylfsApp6}, the required Eq.\eqref{eq:ylfs7b} would be obtained as
\begin{equation}
    \xi^{2n}-\frac{\lambda}{n}\xi^{n+1}-a=0.
    \label{eq:ylfsApp11}
\end{equation}
\appendix

\end{document}